\documentclass{amsproc}

\usepackage[swedish,english]{babel}
\usepackage[utf8]{inputenc}
\usepackage[T1]{fontenc}
\usepackage{amsmath,amssymb,amsfonts}
\usepackage{url}
\usepackage[dvipsnames,usenames]{color}
\usepackage{array,longtable}
\usepackage[pdftex]{graphicx}
\usepackage{tikz}
\usepackage{accents}
\usepackage{tabu}

\newcolumntype{L}{>{\displaystyle}l}
\newcolumntype{C}{>{\displaystyle}c}
\newcolumntype{R}{>{\displaystyle}r}
\newcolumntype{J}{>{$\displaystyle}l<{$}}
\newcolumntype{D}{>{$\displaystyle}c<{$}}
\newcolumntype{S}{>{$\displaystyle}r<{$}}

\renewcommand{\tfrac}{\genfrac{}{}{}1}

\newcommand{\R}{\mathbb R}
\newcommand{\N}{\mathbb N}
\newcommand{\C}{\mathbb C}
\newcommand{\cont}{\mathcal C}
\renewcommand{\Im}{\mathrm{Im}}

\newcommand{\pois}{\mathcal P}
\def\I{\mathfrak{i}}
\newcommand{\diff}{\mathrm{d}}
\renewcommand{\bar}{\overline}

\newcommand{\HN}{Herglotz-Nevanlinna }
\newcommand{\NV}{Nevanlinna }

\newcommand{\ntto}{\:\scriptsize{\xrightarrow{\vee\:}}\:}
\newcommand{\nntto}{\:\scriptsize{\xrightarrow{\wedge\:}}\:}
\newcommand{\cutN}{(\C\setminus\R)^n}
\newcommand{\cut}{\C\setminus\R}

\newcommand{\ie}{\textit{i.e.}\/ } 
\newcommand{\eg}{\textit{e.g.}\/ } 
\newcommand{\cf}{\textit{cf.}\/ } 

\usepackage{bm}
\renewcommand{\vec}{\bm}

\allowdisplaybreaks

\theoremstyle{definition} 
\newtheorem{define}{Definition}[section]
\newtheorem{example}[define]{Example}

\theoremstyle{plain}
\newtheorem{lemma}[define]{Lemma}
\newtheorem{thm}[define]{Theorem}
\newtheorem{prop}[define]{Proposition}
\newtheorem{coro}[define]{Corollary}

\numberwithin{equation}{section}

\begin{document}

\title[Characterization of the symmetric extension]{An analytic  characterization of the symmetric extension of a \HN function in several variables}

\author{Mitja Nedic}
\address{Mitja Nedic, Department of Mathematics and Statistics, University of Helsinki, PO Box 68, FI-00014 Helsinki, Finland, orc-id: 0000-0001-7867-5874}
\curraddr{}
\email{mitja.nedic@helsinki.fi}
\thanks{\textit{Key words.} \HN functions, Cauchy-type functions, symmetric extension.}

\subjclass[2010]{32A36, 32A99.}

\date{2019-11-27}

\begin{abstract}
In this paper, we derive an analytic characterization of the symmetric extension of a \HN function in several variables. Here, the main tools used are the so-called variable non-dependence property and the symmetry formula satisfied by \HN and Cauchy-type functions. We also provide an extension of the Stieltjes inversion formula for Cauchy-type functions.
\end{abstract}

\maketitle

\section{Introduction}\label{sec:intro}

On the upper half-plane $\C^+ := \{z \in \C~|~\Im[z] > 0\}$, the class of holomorphic functions with non-negative imaginary part plays an important role in many areas of analysis and applications. These functions, called \emph{\HN functions}, appear, to name but a few examples, in the theory of Sturm-Liouville operators and their perturbations \cite{Aronszajn1957,AronszajnBrown1970,Donoghue1965,KacKrein1974}, when studying the classical moment problem \cite{Akhiezer1965,Nevanlinna1922,Simon1998}, when deriving physical bounds for passive systems \cite{BernlandEtal2011} or as approximating functions in certain convex optimization problems \cite{IvanenkoETAL2019a,IvanenkoETAL2020}.

A classical integral representation theorem \cite{Cauer1932,Nevanlinna1922} states that any \HN function $h$ can be written, for $z \in \C^+$, as
\begin{equation}
    \label{eq:intRep_1var}
    h(z) = a + bz + \frac{1}{\pi}\int_\R\left(\frac{1}{t-z}-\frac{t}{1+t^2}\right)\diff\mu(t),
\end{equation}
where $a \in \R$, $b \geq 0$ and $\mu$ is a positive Borel measure on $\R$ for which $\int_\R(1+t^2)^{-1}\diff\mu(t) < \infty$. Although this representation is \emph{a prioir} established for $z \in \C^+$, it is well-defined, as an algebraic expression, for any $z \in \cut$. Hence, for a \HN function $h$, we define its \emph{symmetric extension} $h_\mathrm{sym}$ as the right-hand side of representation \eqref{eq:intRep_1var} where we now take $z \in \cut$. It is now an easy consequence of the definitions that a holomorphic function $f\colon\cut\to\C$ equals the symmetric extension of some \HN function if and only if it holds that $\Im[f(z)] \geq 0$ for $z \in \C^+$ and $f(z) = \bar{f(\bar{z})}$ for all $z \in \cut$. In this way, we obtain an analytic characterzation of the symmetric extension.

When considering, instead, functions in the poly-upper half-plane
$$\C^{+n}:= (\C^+)^n = \big\{\vec{z}\in\C^n \,\big |\,\forall j=1,2,\ldots, n:  \Im[z_j]>0 \big\},$$
the analogous situation becomes more involved. \HN functions in several variables, \cf Definition \ref{def:HN_Nvar}, appear \eg when considering operator monotone functions \cite{AglerEtal2012} or with representations of multidimensional passive systems \cite{Vladimirov1979}. Their corresponding integral representation is recalled in detail in Theorem \ref{thm:intRep_Nvar} later on and leads, in an analogous way as in the one-variable case, to the definition of a \emph{symmetric extension}, which is now a holomorphic function on $\cutN$. As such, the main goal of this paper is to give an analytic characterization of symmetric extensions of a \HN function in several variables, \ie we wish to be able to determine when a function $f \colon \cutN \to \C$ is, in fact, equal to the symmetric extension of a \HN function. This is answered by Theorem \ref{thm:symmetric_characterization} and Corollary \ref{coro:symmetric_characterization}. 

The structure of the paper is as follows. After the introduction in Section \ref{sec:intro} we review the different classes of functions that will appear throughout the paper in Section \ref{sec:classes}. Section \ref{sec:symmetry} is then devoted to presenting the main result of the paper as well as some important examples. Finally, Section \ref{sec:Stieltjes_inversion} discusses how the Stieltjes inversion formula can be extended to certain functions on $\cutN$.

\section{Classes of functions in the poly cut-plane}\label{sec:classes}

Throughout this paper, we will primarily consider two classes of holomorphic functions on the poly cut-plane $\cutN$, both of which are intricately connected to a certain kernel function. These objects are defined as follows.

\subsection{The kernel $K_n$ and Cauchy-type functions}\label{subsec:Cauchy_type}

We begin by introducing the kernel $K_n\colon (\C\setminus\R)^n \times \R^n \to \C$ as
\begin{equation}\label{eq:kernel_Nvar}
K_n(\vec{z},\vec{t}) := \I\left(\frac{2}{(2\I)^n}\prod_{\ell=1}^n\left(\frac{1}{t_\ell-z_\ell}-\frac{1}{t_\ell+\I}\right)-\frac{1}{(2\I)^n}\prod_{\ell=1}^n\left(\frac{1}{t_\ell-\I}-\frac{1}{t_\ell+\I}\right)\right).
\end{equation}
If the vector $\vec{z}$ is restricted to $\C^{+n}$, then the kernel $K_n$ is a complex-constant multiple of the Schwartz kernel of $\C^{+n}$ viewed as a tubular domain over the cone $[0,\infty)^n$ \cite[Sec. 12.5]{Vladimirov1979}. 

When $n = 1$, it holds that
$$K_1(z,t) = \frac{1}{t-z} - \frac{t}{1+t^2}.$$
As such, the kernel $K_1$ satisfies, for all $z \in \C\setminus\R$ and all $t \in \R$, the symmetry property
$$K_1(z,t) = \bar{K_1(\bar{z},t)}.$$
When $n \geq 2$, the symmetry satisfied by the kernel becomes more involved and requires the introduction of some additional notation. First, given two numbers $z,w \in \C$, an indexing set $B \subseteq \{1,2,\ldots,n\}$ and an index $j \in \{1,2,\ldots,n\}$, define
$$\psi_B^j(z,w) := \left\{\begin{array}{rcl}
z & ; & j \not\in B, \\
\bar{w} & ; & j \in B.
\end{array}\right.$$
Second, given an indexing set $B \subseteq \{1,2,\ldots,n\}$, define the map $\Psi_B\colon \C^n \times \C^n \to \C^n$ as $\Psi_B(\vec{z},\vec{w}) := \vec{\zeta}$ with $\zeta_j := \psi_B^j(z_j,w_j)$. In other words, the map $\Psi_B$ functions as a way of selectively combining two vectors into one where the set $B$ determines which components of $\vec{z}$ should be replaced by the conjugates of the components of $\vec{w}$. It now holds that
\begin{equation}
    \label{eq:symmetry_Kn}
    K_n(\vec{z},\vec{t}) = \sum_{\substack{B \subseteq \{1,\ldots,n\} \\ B \neq \emptyset}}(-1)^{|B|+1}\bar{K_n(\Psi_B(\I\,\vec{1},\vec{z}),\vec{t})}
\end{equation}
for every $\vec{z} \in (\C\setminus\R)^n$ and every $\vec{t} \in \R^n$ \cite[Prop. 6.1]{LugerNedic2019}.

Using the kernel $K_n$, the largest class of functions that will be considered is the following.

\begin{define}
\label{def:Cauchy_type}
A function $g\colon (\C\setminus\R)^n \to \C$ is called a \emph{Cauchy-type function} if there exists a positive Borel measure $\mu$ on $\R^n$ satisfying the growth condition
\begin{equation}
\label{eq:measure_growth}
\int_{\R^n}\prod_{\ell=1}^n\frac{1}{1+t_\ell^2}\diff\mu(\vec{t}) < \infty
\end{equation}
such that
$$g(\vec{z}) = \frac{1}{\pi^n}\int_{\R^n}K_n(\vec{z},\vec{t})\diff\mu(\vec{t})$$
for every $\vec{z} \in (\C\setminus\R)^n$.
\end{define}

Note that this definition is different from \cite[Def. 3.1]{LugerNedic2019} in that it assumes from the beginning that a Cauchy-type function is defined on $(\C\setminus\R)^n$ and not only on $\C^{+n}$. Furthermore, it would be possible to define an even larger class of functions using the same kernel, but general distributions instead of measures, see \cite[Ex. 7.7]{LugerNedic2019b} for an example. However, this extension will not be considered here. Moreover, Definition \ref{def:Cauchy_type} allows, in principle, for two (or more) different measure to yield the same function $g$, though we will show that this is not the case later in Section \ref{sec:Stieltjes_inversion}.

An immediate consequence of the symmetry formula \eqref{eq:symmetry_Kn} is an analogous symmetry formula for Cauchy-type functions. In particular, it holds, for any Cauchy-type function $g$, that
\begin{equation}
    \label{eq:symmetry_Cauchy_type}
    g(\vec{z}) = \sum_{\substack{B \subseteq \{1,\ldots,n\} \\ B \neq \emptyset}}(-1)^{|B|+1}\bar{g(\Psi_B(\I\,\vec{1},\vec{z}))}
\end{equation}
for every $\vec{z} \in (\C\setminus\R)^n$ and every $\vec{t} \in \R^n$ \cite[Prop. 6.5]{LugerNedic2019}.

The growth of a Cauchy-type function along a coordinate parallel complex line can be described using non-tangential limits. These are taken in so-called \emph{Stoltz domains} and are defined as follows. An \emph{upper Stoltz domain} with centre $0 \in \R$ and angle $\theta \in (0,\frac{\pi}{2}]$ is the set $\{z \in \C^+~|~\theta \leq \arg(z) \leq \pi-\theta\}$ and the symbol $z \ntto \infty$ then denotes the limit $|z| \to \infty$ in any upper Stoltz domain with centre $0$. A \emph{lower Stoltz domain} and the symbol $z \nntto \infty$ are defined analogously. Furthermore, we note that in the literature, slightly different notations are sometimes used to describe these limits. Two examples of Stlotz domains are visualized in Figure \ref{fig:Stoltz_domain} below.

\begin{figure}[!ht]
\centering
\begin{tikzpicture}[scale=1.1]
\fill[fill=black!10!white] (-3,2.179) -- (0,0) -- (3,2.179);
\draw[help lines,->] (-3.2,0) -- (3.2,0) node[above] {$x$};
\draw[help lines,->] (0,-2.5) -- (0,2.5) node[right] {$\I\,y$};
\draw[-] (-3,2.179) -- (0,0) -- (3,2.179);
\draw [dashed,domain=0:36] plot ({0.75*cos(\x)}, {0.75*sin(\x)});
\draw (0.9,0.3) node {$\theta_1$};
\end{tikzpicture}
\begin{tikzpicture}[scale=1.1]
\fill[fill=black!10!white] (-1.7,-2.179) -- (0,0) -- (1.7,-2.179);
\draw[help lines,->] (-1.7,0) -- (1.7,0) node[above] {$x$};
\draw[help lines,->] (0,-2.5) -- (0,2.5) node[right] {$\I\,y$};
\draw[-] (-1.7,-2.179) -- (0,0) -- (1.7,-2.179);
\draw [dashed,domain=0:51] plot ({-0.75*cos(\x)}, {-0.75*sin(\x)});
\draw (-0.9,-0.25) node {$\theta_2$};
\end{tikzpicture}
\caption{An upper Stoltz domain with centre $0$ and angle $\theta_1$ (left) and a lower Stoltz domain with centre $0$ and angle $\theta_2$ (right).}
\label{fig:Stoltz_domain}
\end{figure}
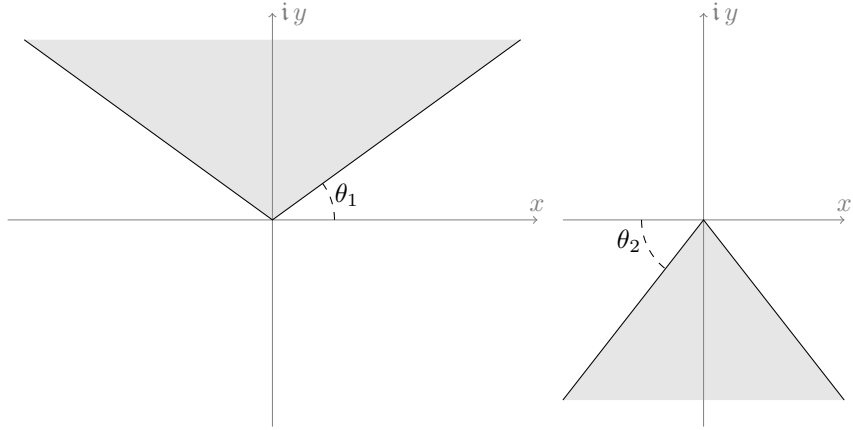

For any Cauchy-type function $g$ it now holds, for any $\vec{z} \in \cutN$ and any $j \in \{1,\ldots,n\}$, that
$$\lim\limits_{z_j \ntto \infty}\frac{g(\vec{z})}{z_j} = \lim\limits_{z_j \nntto \infty}\frac{g(\vec{z})}{z_j} = 0,$$
see \cite[Lem. 3.2]{LugerNedic2019}.

\subsection{Herglotz-Nevanlina functions}\label{subsec:HN}

These functions are defined as follows, \cf \cite{LugerNedic2017,LugerNedic2019,Vladimirov1969,Vladimirov1979}.

\begin{define}\label{def:HN_Nvar}
A holomorphic function $h\colon \C^{+n} \to \C$ is called a \emph{\HN function} if it is holomorphic with non-negative imaginary part.
\end{define}

In contrast to the definition of Cauchy-type function, the above definition is analytic in nature, \ie it describes the function class in terms of conditions on the function itself. In order to be able to relate it to the kernel $K_n$, we introduce, given ambient numbers $z \in \C\setminus\R$ and $t \in \R$, the expressions
$$\begin{array}{RCL}
N_{-1}(z,t) & := & \frac{1}{2\,\I}\left(\frac{1}{t - z} - \frac{1}{t - \I}\right), \\[0.35cm]
N_{0}(z,t) & := & \frac{1}{2\,\I}\left(\frac{1}{t - \I} - \frac{1}{t_j + \I}\right), \\[0.35cm]
N_{1}(z,t) & := & \frac{1}{2\,\I}\left(\frac{1}{t + \I} - \frac{1}{t - \bar{z}}\right).
\end{array}$$
Note that $N_0$ is independent of $z \in \C\setminus\R$ and $N_0(z,t) \in \R$ while
$$\bar{N_{-1}(z,t)} = N_{1}(z,t)$$
for all $z \in \C\setminus\R$ and $t \in \R$. Using these expression, one may give an integral representation formula for \HN functions involving the kernel $K_n$ \cite[Thm. 4.1]{LugerNedic2019}.

\begin{thm}\label{thm:intRep_Nvar}
A function $h\colon \C^{+n} \to \C$ is a \HN function if and only if $h$ can be written as
\begin{equation}\label{eq:intRep_Nvar}
h(\vec{z}) = a + \sum_{j=1}^nb_j\,z_j + \frac{1}{\pi^n}\int_{\R^n}K_n(\vec{z},\vec{t})\diff\mu(\vec{t}),
\end{equation}
where $a \in \R$, $\vec{b} \in [0,\infty)^n$, the kernel $K_n$ is as before and $\mu$ is a positive Borel measure on $\R^n$ satisfying the growth condition \eqref{eq:measure_growth} and the \NV condition
\begin{equation}
\label{eq:measure_Nevan}
\sum_{\substack{\vec{\rho} \in \{-1,0,1\}^n \\ -1 \in \vec{\rho} \wedge 1 \in \vec{\rho}}}\int_{\R^n}\prod_{j=1}^nN_{\rho_j}(z_j,t_j)\diff\mu(\vec{t}) = 0
\end{equation}
for all $\vec{z} \in \C^{+n}$. Furthermore, for a given function $h$, the triple of representing parameters $(a,\vec{b},\mu)$ is unique.
\end{thm}

The integral representation in formula \eqref{eq:intRep_Nvar} is well-defined for any $\vec{z} \in (\C\setminus\R)^n$, which may be used to extend any Herglotz-Nevanlinna function $h$ from $\C^{+n}$ to $(\C\setminus\R)^n$. This extension is called the \emph{symmetric extension} of the function $h$ and is denoted as $h_\mathrm{sym}$. The symmetric extension of a \HN function $h$ is different from its possible analytic extension as soon as $\mu \neq 0$ \cite[Prop. 6.10]{LugerNedic2019} and satisfies the following variable-dependence property \cite[Prop. 6.9]{LugerNedic2019}.

\begin{prop}\label{prop:variable_dependence}
    Let $n \geq 2$ and let $h_\mathrm{sym}$ be the symmetric extension of a Herglotz-Nevanlinna function $h$ in $n$ variables for which $\vec{b} = \vec{0}$. Let $\vec{z} \in (\C\setminus\R)^n$ be such that $z_j \in \C^-$ for some index $j \in \{1,2,\ldots,n\}$. Then, the value $h_\mathrm{sym}(\vec{z})$ does not depend on the components of $\vec{z}$ that lie in $\C^+$.
\end{prop}

Furthermore, if $h$ is a \HN function for which $\vec{b} = \vec{0}$, then its symmetric extension $h_\mathrm{sym}$ will satisfy the symmetry formula
\begin{equation}
    \label{eq:symmetry_HN}
    h_\mathrm{sym}(\vec{z}) = \sum_{\substack{B \subseteq \{1,\ldots,n\} \\ B \neq \emptyset}}(-1)^{|B|+1}\bar{h_\mathrm{sym}(\Psi_B(\I\,\vec{1},\vec{z}))},
\end{equation}
where $\vec{z} \in (\C\setminus\R)^n$ \cite[Prop. 6.7]{LugerNedic2019}. When $n = 1$, it is not necessary to assume that $b = 0$ for formula \eqref{eq:symmetry_HN} to hold. However, when $n > 1$, this is required.

The representing vector $\vec{b}$ describes the growth of the function $h$ along coordinate parallel complex lines in $\C^{+n}$. More precisely, we recall from \cite[Cor. 4.6(iv)]{LugerNedic2019} that, for any $j \in \{1,\ldots,n\}$, we have
\begin{equation}
    \label{eq:b_parameters}
    b_j  = \lim\limits_{z_j \ntto \infty}\frac{h(\vec{z})}{z_j}.
\end{equation}
In particular, the above limit is independent of the entries of the vector $\vec{z}$ at the non-$j$-th positions. This result carries over to the symmetric extension, for which it holds, for any $j \in \{1,\ldots,n\}$, that
$$b_j = \lim\limits_{z_j \ntto \infty}\frac{h_\mathrm{sym}(\vec{z})}{z_j} = \lim\limits_{z_j \nntto \infty}\frac{h_\mathrm{sym}(\vec{z})}{z_j}.$$

Every \HN function that is represented by a data-triple of the form $(0,\vec{0},\mu)$ in the sense of Theorem \ref{thm:intRep_Nvar} is also a Cauchy-type function. The converse, \ie that every Cauchy-type function equals a \HN function represented by a data-triple of the form $(0,\vec{0},\mu)$, is true only when $n = 1$. This is due to the fact that when $n = 1$, the Nevanlinna condition \eqref{eq:measure_Nevan} becomes emptily fulfilled by every positive Borel measures $\mu$ satisfying the growth condition \eqref{eq:measure_growth}.

\section{Symmetry and variable non-dependence}
\label{sec:symmetry}

We begin by recalling that the symmetric extension of a \HN function $h$ in one variable is uniquely determined by its values in $\C^+$. Indeed, when $n = 1$, the symmetry formula \eqref{eq:symmetry_HN} takes the form $h_\mathrm{sym}(z) = \bar{h_\mathrm{sym}(\bar{z})}$, providing a way to recover the values of the function in $\C^-$ using only the values of the function in $\C^+$.

For functions of several variables, the appropriate analogue involves the following definition.

\begin{define}
\label{def:variable_non_dependance}
A function $f \colon \cutN \to \C$ is said to satisfy the \emph{variable non-dependence property} if for every vector $\vec{z} \in \cutN$ such that $z_j \in \C^-$ for some index $j \in \{1,2,\ldots,n\}$ the value $f(\vec{z})$ does not depend on the components of $\vec{z}$ that lie in $\C^+$.
\end{define}

By Proposition \ref{prop:variable_dependence}, the symmetric extension of a \HN function satisfies the variable non-dependence property \ref{def:variable_non_dependance} if $\vec{b} = \vec{0}$. In particular, the symmetric extension of any \HN function that is also a Cauchy-type function will always satisfy the variable non-dependence property \ref{def:variable_non_dependance}. However, a general Cauchy-type function need not satisfy it, as shown by the function $f_2$ in Example \ref{ex:independent_conditions} later on.

We may now describe the precise circumstances under which we can recover the values of a function defined on $\cutN$ purely in terms of its values in $\C^{+n}$.

\begin{prop}\label{prop:symmetry_uniqueness}
Let $f \colon \cutN \to \C$ be a holomorphic function satisfying the symmetry formula \eqref{eq:symmetry_HN} and the variable non-dependence property \ref{def:variable_non_dependance}. Then, the values of the function $f$ on $\cutN$ are uniquely determined by its values in $\C^{+n}$.
\end{prop}

\proof
Using the symmetry formula \eqref{eq:symmetry_HN}, let us investigate the values of the function $f$ in a connected component of $(\C\setminus\R)^n$ where at least one of the coordinates has a negative sign of the imaginary part, \ie we are investigating a connected component $X \subseteq (\C\setminus\R)^n$ where exist at lest one index $j \in \{1,\ldots,n\}$ such that the $j$-th coordinate lies in $\C^-$. For any such chosen connected component $X$, let $B' \subseteq \{1,\ldots,n\}$ be the set of those indices for which the corresponding variables lie in $\C^-$. In particular, $1 \leq |B'| \leq n$. For $\vec{z} \in X$, it holds, by the symmetry formula \eqref{eq:symmetry_HN}, that
\begin{multline*}
    f(\vec{z}) = \sum_{\substack{B \subseteq \{1,\ldots,n\} \\ B \neq \emptyset}}(-1)^{|B|+1}\bar{f(\Psi_B(\I\,\vec{1},\vec{z}))} \\
    = \sum_{\substack{B \subseteq \{1,\ldots,n\} \\ B \neq \emptyset \wedge B \subseteq B'}}(-1)^{|B|+1}\bar{f(\Psi_B(\I\,\vec{1},\vec{z}))} + \sum_{\substack{B \subseteq \{1,\ldots,n\} \\ B \not\subseteq B'}}(-1)^{|B|+1}\bar{f(\Psi_B(\I\,\vec{1},\vec{z}))}.
\end{multline*}
Due to the definition of the set $B'$, it holds that $\Psi_{B}(\I\,\vec{1},\vec{z}) \in \C^{+n}$ for any $\vec{z} \in X$ and any indexing set $B \subseteq B'$. Furthermore, by the variable non-dependence property \ref{def:variable_non_dependance}, it holds that
$$\Psi_{B}(\I\,\vec{1},\vec{z}) = \Psi_{B\setminus B'}(\I\,\vec{1},\vec{z})$$
for any $\vec{z} \in X$ and any indexing set $B$ where $B \not\subseteq B'$. Hence,
$$f(\vec{z}) = \sum_{\substack{B \subseteq \{1,\ldots,n\} \\ B \neq \emptyset \wedge B \subseteq B'}}(-1)^{|B|+1}\bar{f(\Psi_B(\I\,\vec{1},\vec{z}))} + \sum_{\substack{B \subseteq \{1,\ldots,n\} \\ B \not\subseteq B'}}(-1)^{|B|+1}\bar{f(\Psi_{B\setminus B'}(\I\,\vec{1},\vec{z}))}.$$

We now claim that the second sum is always equal to zero. Indeed, if $|B'| = n$, there is nothing left to prove. Otherwise, we may assume that $|B'| < n$, where we claim that there is a way to "pair up" the indexing sets in the second sum in such a way that the two sets in each pair only differ by one element in $B'$. We construct this pairing in the following way. Let $j_1$ be the smallest index in $B'$. Then, exactly half of the sets $B \subseteq \{1,\ldots,n\}$ that are not subsets of $B'$ contain the index $j_1$ and exactly half of them do not contain the index $j_1$. This follows from the general observation that exactly half of the subsets of a given set contain a specific element of the set. An indexing set $B_1$ is then paired with the indexing set $B_1 \cup \{j_1\}$. In this case, $(B_1 \cup \{j_1\})\setminus B' = B_1 \setminus B'$ and
\begin{multline*}
    (-1)^{|B_1|+1}\bar{f(\Psi_{B_1\setminus B'}(\I\,\vec{1},\vec{z}))} + (-1)^{|B_1 \cup \{j_1\}|+1}\bar{f(\Psi_{(B_1 \cup \{j_1\})\setminus B'}(\I\,\vec{1},\vec{z}))} \\
    = (-1)^{|B_1|+1}\bar{f(\Psi_{B_1\setminus B'}(\I\,\vec{1},\vec{z}))} - (-1)^{|B_1|+1}\bar{f(\Psi_{B_1 \setminus B'}(\I\,\vec{1},\vec{z}))} = 0,
\end{multline*}
yielding the desired result.
\endproof

Using Proposition \ref{prop:symmetry_uniqueness}, we may now given an analytic characterization of the symmetric extension of a \HN function.

\begin{thm}
\label{thm:symmetric_characterization}
Let $f \colon \cutN \to \C$ be a holomorphic function such that
$$\lim\limits_{z_j \ntto \infty}\frac{f(\vec{z})}{z_j} = \lim\limits_{z_j \nntto \infty}\frac{f(\vec{z})}{z_j} = 0$$
for all indices $j \in \{1,\ldots,n\}$. Then $f = h_\mathrm{sym}$ for some \HN function $h$ if and only if
\begin{itemize}
    \item[(i)]{it holds that $\Im[f(\vec{z})] \geq 0$ for all $\vec{z} \in \C^{+n}$,}
    \item[(ii)]{the function $f$ satisfies the symmetry formula \eqref{eq:symmetry_HN},}
    \item[(iii)]{the function $f$ satisfies the variable non-dependence property \ref{def:variable_non_dependance}.}
\end{itemize}
\end{thm}

\proof
If $f = h_\mathrm{sym}$ for some \HN function $h$, then this function must have $\vec{b} = \vec{0}$ due to the assumption on the growth of $f$. Then, properties (i) -- (iii) are satisfied by the previously known results discussed in Section \ref{subsec:HN}. Conversely, if we are given the function $f$ satisfies the properties (i) -- (iii), we construct a \HN function out the function $f$ by setting
$$h := f|_{\C^{+n}}.$$
This function $h$ may then be symmetrically extended to $\cutN$. However, $f$ and $h_\mathrm{sym}$ are now two holomorphic functions on $\cutN$ satisfying the symmetry formula \eqref{eq:symmetry_HN} and the variable non-dependence property \ref{def:variable_non_dependance} which, furthermore, agree on $\C^{+n}$. Therefore, by Proposition \ref{prop:symmetry_uniqueness}, they agree everywhere on $\cutN$, as desired.
\endproof

The assumption on the growth of the function $f$ may be slightly weakened, but, to compensate, conditions (ii) and (iii) need to be slightly modified. 

\begin{coro}
\label{coro:symmetric_characterization}
Let $f \colon \cutN \to \C$ be a holomorphic function such that
$$\lim\limits_{z_j \ntto \infty}\frac{f(\vec{z})}{z_j} = \lim\limits_{z_j \nntto \infty}\frac{f(\vec{z})}{z_j} = d_j \geq 0$$
for all indices $j \in \{1,\ldots,n\}$. In particular, for a fixed $j \in \{1,\ldots,n\}$, the above limits are assumed to be independent of the values of the vector $\vec{z} \in \cutN$ at the non-$j$-th positions. Then $f = h_\mathrm{sym}$ for some \HN function $h$ if and only if
\begin{itemize}
    \item[(i)]{it holds that $\Im[f(\vec{z})] \geq 0$ for all $\vec{z} \in \C^{+n}$,}
    \item[(ii')]{the function $\vec{z} \mapsto f(\vec{z}) - \sum_{j=1}^nd_jz_j$ satisfies the symmetry formula \eqref{eq:symmetry_HN},}
    \item[(iii')]{the function $\vec{z} \mapsto f(\vec{z}) - \sum_{j=1}^nd_jz_j$ satisfies the variable non-dependence property \ref{def:variable_non_dependance}.}
\end{itemize}
\end{coro}

The three conditions on the function $f$ in Theorem \ref{thm:symmetric_characterization} are independent of each-other. To verify this, consider the following functions on $\cutN$.

\begin{example}
\label{ex:independent_conditions}
Table \ref{tab:dependence} presents eight explicit functions defined on $\cutN$ and Table \ref{tab:dependence_part2} summarizes which conditions of Theorem \ref{thm:symmetric_characterization} are fulfilled by which function. Note also that all eight functions satisfy the assumption on the growth of the function from Theorem \ref{thm:symmetric_characterization}. The functions are constructed as follows.

The functions $f_0$ is defined to equal: a negative imaginary constant on $\C^+ \times \C^+$, breaking condition (i); a function depending only on the second variable on $\C^- \times \C^+$, breaking condition (iii); and identically zero in the remaining connected components of $\cutN$, ensuring that condition (ii) is not satisfied. The function $f_1$ is obtained form $f_0$ by changing the definition on $\C^+ \times \C^+$ to a positive imaginary constant, thereby satisfying condition (i), but still neither (ii) nor (iii).

The function $f_2$ is the Cauchy-type function given by a measure $\mu_2$ on $\R^2$ defined on Borel subsets $U \subseteq \R^2$ as
$$\mu_2(U) := \pi\int_{\R}\chi_U(t,t)\diff t,$$
where $\chi$ denotes the characteristic function of a set. This measure obviously satisfies the growth condition \eqref{eq:measure_growth} and it does not satisfy the Nevanlinna condition \eqref{eq:measure_Nevan} as it supported on the diagonal in $\R^2$ - an impossibility for Nevanlinna measures as shown in \cite[Ex. 3.14]{LugerNedic2020}. This function does not satisfy condition (i) as, for example, $f_2(4\,\I,4\,\I) = -\frac{\I}{10}$. As a Cauchy-type function, it is guaranteed to satisfy condition (ii). It also clearly does not satisfy condition (iii) as the values in \eg $\C^+ \times \C^-$ depend explicitly on both variables. Note now that while the function $f_2$ takes values with negative imaginary part in $\C^+ \times \C^+$, its imaginary part is bounded form below. Indeed, the functions $z_1 \mapsto -\frac{1}{\I+z_1}$ and $z_2 \mapsto -\frac{1}{\I+z_2}$ are \HN functions of one variable, implying that $\Im[f_2(z_1,z_2)] \geq -\frac{1}{2}$ for all $(z_1,z_2) \in \C^+ \times \C^+$. Hence, the function $f_4$ is obtained by adding to the function $f_2$ the symmetric extension of the \HN function $(z_1,z_2) \mapsto 5\,\I$ (represented by the measure $5\lambda_{\R^2}$). This new function now satisfies condition (i) in addition to (ii), while clearly still not satisfying condition (iii). Note that the function
$$f_4|_{\C^+ \times \C^+}(z_1,z_2) = \frac{9\,\I}{2} - \frac{1}{\I+z_1} - \frac{1}{\I+z_2}$$
as a \HN function is not represented by the measure $\mu_2 + 5\lambda_{\R^2}$ in the sense of Theorem \ref{thm:intRep_Nvar}, but rather by the measure
$$\tfrac{9}{2}\lambda_{\R^2} + (\tau \mapsto (1+\tau^2)^{-1})\lambda_\R \otimes \lambda_\R + \lambda_\R \otimes (\tau \mapsto (1+\tau^2)^{-1})\lambda_\R.$$

The function $f_3$ is defined as zero on all the connected components of $(\cut)^2$ other than $\C^+\times\C^+$ to ensure that it satisfies condition (iii), while setting the function equal to a negative imaginary constant in $\C^+ \times \C^+$ ensures that it satisfies neither condition (i) nor (ii). Changing this definition to a positive imaginary constant in $\C^+ \times \C^+$ gives the function $f_5$ which satisfies condition (i) and (iii), but not (ii).

The function $f_7$ is simply taken as the symmetric extension of a \HN function, thereby satisfying all three properties automatically. Finally, the function $f_6$ is chosen as $f_6 := -f_7$, satisfying conditions (ii) and (iii), but not (i). \hfill$\lozenge$
\end{example}

\begin{table}[!ht]
    \centering
    \begin{tabular}{c|cccc}
~ & $\C^+ \times \C^+$ & $\C^- \times \C^+$ & $\C^+ \times \C^-$ & $\C^- \times \C^-$  \\
\hline
$f_0$ & $-\I$ & $\frac{1}{z_2}$ & $0$ & $0$  \\[0.2cm] 
$f_1$ & $\I$ & $\frac{1}{z_2}$ & $0$ & $0$  \\[0.2cm] 
$f_2$ & $-\frac{\I}{2} - \frac{1}{\I+z_1} - \frac{1}{\I+z_2}$ & $-\frac{\I}{2} + \frac{1}{z_2-z_1} - \frac{1}{\I+z_2}$ & $-\frac{\I}{2} - \frac{1}{\I+z_1} + \frac{1}{z_1-z_2}$ & $-\frac{\I}{2}$  \\[0.2cm] 
$f_3$ & $-\I$ & $0$ & $0$ & $0$  \\[0.2cm] 
$f_4$ & $\frac{9\,\I}{2} - \frac{1}{\I+z_1} - \frac{1}{\I+z_2}$ & $-\frac{11\,\I}{2} + \frac{1}{z_2-z_1} - \frac{1}{\I+z_2}$ & $-\frac{11\,\I}{2} - \frac{1}{\I+z_1} + \frac{1}{z_1-z_2}$ & $-\frac{11\,\I}{2}$  \\[0.2cm] 
$f_5$ & $\I$ & $0$ & $0$ & $0$  \\[0.2cm] 
$f_6$ & $-\I$ & $\I$ & $\I$ & $\I$  \\[0.2cm] 
$f_7$ & $\I$ & $-\I$ & $-\I$ & $-\I$  
    \end{tabular}
    \caption{Examples of eight functions defined on $(\cut)^2$.}
    \label{tab:dependence}
\end{table}

\begin{table}[!ht]
    \centering
    \begin{tabular}{c|ccc}
~ & (i) & (ii) & (iii) \\
\hline
$f_0$ & $\times$ & $\times$ & $\times$ \\ 
$f_1$ & \checkmark & $\times$ & $\times$ \\ 
$f_2$ & $\times$ & \checkmark & $\times$ \\ 
$f_3$ & $\times$ & $\times$ & \checkmark \\ 
$f_4$ & \checkmark & \checkmark & $\times$ \\ 
$f_5$ & \checkmark & $\times$ & \checkmark \\ 
$f_6$ & $\times$ & \checkmark & \checkmark \\ 
$f_7$ & \checkmark & \checkmark & \checkmark 
    \end{tabular}
    \caption{The relation of the eight functions from Table \ref{tab:dependence} to the three conditions from Theorem \ref{thm:symmetric_characterization}.}
    \label{tab:dependence_part2}
\end{table}

\section{The Stieltjes inversion formula for Cauchy-type functions}
\label{sec:Stieltjes_inversion}

For \HN functions, the Stieltjes inversion formula describes how to reconstruct the representing measure $\mu$ of a \HN function $h$ from the values of the imaginary part of the function in $\C^{+n}$. More precisely, it holds that
$$\int_{\R^n}\varphi(\vec{t})\diff\mu(\vec{t}) = \lim\limits_{\vec{y} \to \vec{0}^+}\int_{\R^n}\varphi(\vec{x})\Im[h(\vec{x} + \I\,\vec{y})]\diff\vec{x}$$
for all $\cont^1$-functions $\varphi\colon\R^n \to \R$ for which there exists a constant $D \geq 0$ such that $|\varphi(\vec{x})| \leq D\prod_{j=1}^n(1+x_j^2)^{-1}$ for all $\vec{x} \in \R^n$, see \eg \cite{KacKrein1974} or \cite[Lem. 4.1]{BernlandEtal2011} for the case $n=1$ and \cite[Cor. 4.6(viii)]{LugerNedic2019} for the general case.

As noted in Section \ref{subsec:HN}, Cauchy-type functions are a subclass of \HN functions when $n=1$ and, hence, one only need the values of (the imaginary part of) a Cauchy-type function in $\C^+$ to reconstruct its measure. However, in Example \ref{ex:independent_conditions}, we have seen two different positive Borel measures on $\R^2$ for which the corresponding Cauchy-type functions agree on $\C^{+2}$, but not on the remaining connected components of $(\cut)^2$. 

The crucial role in the proof of the Stieltjes inversion formula is held by the Poisson kernel of $\C^{+n}$, which, we recall, is defined for $\vec{z} \in \C^{+n}$ and $\vec{t} \in \R^n$ as
$$\pois_n(\vec{z},\vec{t}) := \prod_{j=1}^n\frac{\Im[z_j]}{|t_j-z_j|^2}.$$
Note that $\pois_n(\vec{z},\vec{t}) > 0$ for every $\vec{z} \in \C^{+n}$ and $\vec{t} \in \R^n$. The imaginary part of the kernel $K_n$ is equal to the Poisson kernel $\pois_n$ plus a remainder term which can be expressed in terms of the $N_j$-factors \cite[Prop. 3.3]{LugerNedic2019} and the integral of the remainder with respect to any Nevanlinna measure is zero. 

The following lemma now shows how one can recover the value of the Poison kernel $\pois_n$ at some point $\vec{z} \in \C^{+n}$ (and $\vec{t} \in \R^n$) using the values of kernel $K_n$ form all of the connected components of the poly cut-plane $(\C\setminus\R)^n$.

\begin{lemma}\label{lem:Stieltjes_inversion}
Let $n \in \N$, $\vec{z} \in \C^{+n}$ and $\vec{t} \in \R^n$. Then, it holds that
$$2\,\I\,\pois_n(\vec{z},\vec{t}) = \sum_{B \subseteq \{1,\ldots,n\}}(-1)^{|B|}K_n(\Psi_B(\vec{z},\vec{z}),\vec{t}),$$
where $\Psi_B$ is the selective conjugation map from Section \ref{subsec:Cauchy_type}.
\end{lemma}

\proof
The proof is done by induction on the dimension $n$. If $n = 1$, then
\begin{multline*}
    \sum_{B \subseteq \{1\}}(-1)^{|B|}K_1(\Psi_B(z,z),t) = K_1(\Psi_\emptyset(z,z),t) + (-1)K_1(\Psi_{\{1\}}(z,z),t) \\
    = K_1(z,t) - K_1(\bar{z},t) = 2\,\I\,\Im[K_1(z,t)] = 2\,\I\,\pois_1(z,t),
\end{multline*}
as desired.

Assume now that the statement of the lemma holds for all $n=1,2,\ldots,N-1$ for some $N \in \N$. For $n = N$, take $\vec{z} \in \C^{+N}$ and $\vec{t} \in \R^N$ and let $\vec{z}'$ and $\vec{t}'$ denote the same vectors with the last component removed, \ie $\vec{z}' := (z_1,\ldots,z_{N-1})$ and $\vec{t}' := (t_1,\ldots,t_{N-1})$. Furthermore, denote
$$A(z,t) := \frac{1}{2\,\I}\left(\frac{1}{t-z}-\frac{1}{t+\I}\right).$$
Then, we then calculate that
\begin{longtable}{SDJ}
\multicolumn{3}{J}{\sum_{B \subseteq \{1,\ldots,N\}}(-1)^{|B|}K_{N}(\Psi_B(\vec{z},\vec{z}),\vec{t})} \\[0.6cm]
~ & = & \sum_{\substack{B \subseteq \{1,\ldots,N\} \\ N \not\in B}}(-1)^{|B|}K_{N}(\Psi_B(\vec{z},\vec{z}),\vec{t}) + \sum_{\substack{B \subseteq \{1,\ldots,N\} \\ N \in B}}(-1)^{|B|}K_{N}(\Psi_B(\vec{z},\vec{z}),\vec{t}) \\[0.6cm]
~ & = & \sum_{B' \subseteq \{1,\ldots,N-1\}}(-1)^{|B'|}\bigg[\I\bigg(2\prod_{j=1}^{N-1}A(\psi_{B'}^j(z_j,z_j),t_j)\cdot A(z_N,t_N) - \prod_{j=1}^NA(\I,t_j)\bigg)\bigg] \\[0.6cm]
~ & ~ & +\sum_{B' \subseteq \{1,\ldots,N-1\}}(-1)^{|B'|+1}\bigg[\I\bigg(2\prod_{j=1}^{N-1}A(\psi_{B'}^j(z_j,z_j),t_j)\cdot A(\bar{z}_N,t_N) - \prod_{j=1}^NA(\I,t_j)\bigg)\bigg] \\[0.6cm]
~ & = & 2\,\I\,A(z_N,t_N)\sum_{B' \subseteq \{1,\ldots,N-1\}}(-1)^{|B'|}K_{N-1}(\Psi_{B'}(\vec{z}',\vec{z}'),\vec{t}') \\[0.6cm]
~ & ~ & + \,\I\,\prod_{j=1}^{N-1}A(\I,t_j)\cdot(A(z_N,t_N)-A(\I,t_N))\cdot\overbrace{\sum_{B' \subseteq \{1,\ldots,N-1\}}(-1)^{|B'|}}^{=0} \\[0.6cm]
~ & ~ & -2\,\I\,A(\bar{z}_N,t_N)\sum_{B' \subseteq \{1,\ldots,N-1\}}(-1)^{|B'|}K_{N-1}(\Psi_{B'}(\vec{z}',\vec{z}'),\vec{t}') \\[0.6cm]
~ & ~ & - \,\I\,\prod_{j=1}^{N-1}A(\I,t_j)\cdot(A(\bar{z}_N,t_N)-A(\I,t_N))\cdot\overbrace{\sum_{B' \subseteq \{1,\ldots,N-1\}}(-1)^{|B'|}}^{=0} \\[0.6cm]
~ & = & 2\,\I\,A(z_N,t_N)\,\pois_{N-1}(\vec{z}',\vec{t}') - 2\,\I\,A(\bar{z}_N,t_N)\,\pois_{N-1}(\vec{z}',\vec{t}') = 2\,\I\,\pois_N(\vec{z},\vec{t}),
\end{longtable}
\noindent finishing the proof.
\endproof

The Stieltjes inversion for Cauchy-type functions is, thus, the following.

\begin{thm}\label{thm:Cauchy_type_Stieltjes}
Let $g$ be a Cauchy-type function given by a measure $\mu$. Then, it holds that
\begin{multline}\label{eq:Cauchy_type_Stieltjes}
    \int_{\R^n}\varphi(\vec{t})\diff\mu(\vec{t}) \\ = \lim\limits_{\vec{y} \to \vec{0}^+}\frac{1}{2\,\I}\int_{\R^n}\varphi(\vec{x})\bigg[\sum_{B \subseteq \{1,\ldots,n\}}(-1)^{|B|}g(\Psi_B(\vec{x}+\I\,\vec{y},\vec{x}+\I\,\vec{y}))\bigg]\diff\vec{x}
\end{multline}
for all $\cont^1$-functions $\varphi\colon\R^n \to \R$ for which there exists a constant $D \geq 0$ such that $|\varphi(\vec{x})| \leq D\prod_{j=1}^n(1+x_j^2)^{-1}$ for all $\vec{x} \in \R^n$.
\end{thm}

\proof
By the definition of Cauchy-type functions and Lemma \ref{lem:Stieltjes_inversion}, it holds that
\begin{multline*}
    \sum_{B \subseteq \{1,\ldots,n\}}(-1)^{|B|}g(\Psi_B(\vec{x}+\I\,\vec{y},\vec{x}+\I\,\vec{y})) \\
    = \frac{1}{\pi^n}\int_{\R^n}\bigg[\sum_{B \subseteq \{1,\ldots,n\}}(-1)^{|B|}K_n(\Psi_B(\vec{x}+\I\,\vec{y},\vec{x}+\I\,\vec{y}),\vec{t})\bigg]\diff\mu(\vec{t}) \\
    = \frac{2\,\I}{\pi^n}\int_{\R^n}\pois_{n}(\vec{x} + \I\,\vec{y},\vec{t})\diff\mu(\vec{t}).
\end{multline*}
Hence,
\begin{multline*}
    \lim\limits_{\vec{y} \to \vec{0}^+}\frac{1}{2\,\I}\int_{\R^n}\varphi(\vec{x})\bigg[\sum_{B \subseteq \{1,\ldots,n\}}(-1)^{|B|}g(\Psi_B(\vec{x}+\I\,\vec{y},\vec{x}+\I\,\vec{y}))\bigg]\diff\vec{x} \\
    = \lim\limits_{\vec{y} \to \vec{0}^+}\frac{1}{\pi^n}\int_{\R^n}\varphi(\vec{x})\bigg(\int_{\R^n}\pois_{n}(\vec{x} + \I\,\vec{y},\vec{t})\diff\mu(\vec{t})\bigg)\diff\vec{x} \\
    = \lim\limits_{\vec{y} \to \vec{0}^+}\frac{1}{\pi^n}\int_{\R^n}\bigg(\int_{\R^n}\varphi(\vec{x})\pois_{n}(\vec{x} + \I\,\vec{y},\vec{t})\diff\vec{x}\bigg)\diff\mu(\vec{t}),
\end{multline*}
where the assumptions on the function $\varphi$ and condition \eqref{eq:measure_growth} for $\mu$ justify the use of Fubini's theorem to change the order of integration. The same assumptions permit for Lebesgue's dominated convergence to be used, allowing us to take the limit as $\vec{y} \to \vec{0}^+$ before integrating with respect to the measure $\mu$. Noting that, by \eg \cite[pg. 111]{Koosis1998},
$$\lim\limits_{\vec{y} \to \vec{0}^+}\int_{\R^n}\varphi(\vec{x})\pois_{n}(\vec{x} + \I\,\vec{y},\vec{t})\diff\vec{x} = \pi^n\varphi(\vec{t})$$
finishes the proof.
\endproof

As an immediate corollary of the previous theorem, we may now establish that the correspondence between a Cauchy-type function and its defining measure $\mu$ is, indeed, a bijection.

\begin{coro}
\label{coro:Cauchy_type_uniqueness}
Let $\mu_1,\mu_2$ be two positive Borel measures on $\R^n$ satisfying the growth condition \eqref{eq:measure_growth}. Then,
$$\int_{\R^n}K_n(\vec{z},\vec{t})\diff\mu_1(\vec{t}) = \int_{\R^n}K_n(\vec{z},\vec{t})\diff\mu_2(\vec{t})$$
for all $\vec{z} \in (\C\setminus\R)^n$ if and only if $\mu_1 \equiv \mu_2$.
\end{coro}

\section*{Acknowledgments}

The author would like to thank Dale Frymark for many enthusiastic discussions on the subject.

\bibliographystyle{amsplain}
\bibliography{MitjaNedic_symmetric_extension}

\providecommand{\bysame}{\leavevmode\hbox to3em{\hrulefill}\thinspace}
\providecommand{\MR}{\relax\ifhmode\unskip\space\fi MR }
\providecommand{\MRhref}[2]{%
  \href{http://www.ams.org/mathscinet-getitem?mr=#1}{#2}
}
\providecommand{\href}[2]{#2}
\begin{thebibliography}{10}

\bibitem{AglerEtal2012}
J.~Agler, J.~E. McCarthy, and N.~J. Young, \emph{Operator monotone functions
  and {L}\"owner functions of several variables}, Ann. of Math. (2)
  \textbf{176} (2012), no.~3, 1783--1826.

\bibitem{Akhiezer1965}
N.~I. Akhiezer, \emph{The classical moment problem and some related questions
  in analysis}, Translated by N. Kemmer, Hafner Publishing Co., New York, 1965.

\bibitem{Aronszajn1957}
N.~Aronszajn, \emph{On a problem of {W}eyl in the theory of singular
  {S}turm-{L}iouville equations}, Amer. J. Math. \textbf{79} (1957), 597--610.

\bibitem{AronszajnBrown1970}
N.~Aronszajn and R.~D. Brown, \emph{Finite-dimensional perturbations of
  spectral problems and variational approximation methods for eigenvalue
  problems. {I}: {F}inite-dimensional perturbations}, Studia Math. \textbf{36}
  (1970), 1--76.

\bibitem{BernlandEtal2011}
A.~Bernland, A.~Luger, and M.~Gustafsson, \emph{Sum rules and constraints on
  passive systems}, J. Phys. A: Math. Theor. \textbf{44} (2011), no.~14,
  145205.

\bibitem{Cauer1932}
W.~Cauer, \emph{The {P}oisson integral for functions with positive real part},
  Bull. Amer. Math. Soc. \textbf{38} (1932), no.~10, 713--717.

\bibitem{Donoghue1965}
W.~F. Donoghue, Jr., \emph{On the perturbation of spectra}, Comm. Pure Appl.
  Math. \textbf{18} (1965), 559--579.

\bibitem{IvanenkoETAL2019a}
Y.~Ivanenko, M.~Gustafsson, B.~L.~G. Jonsson, A.~Luger, B.~Nilsson, S.~Nordebo,
  and J.~Toft, \emph{Passive approximation and optimization using {B}-splines},
  SIAM Journal on Applied Mathematics \textbf{79} (2019), no.~1, 436--458.

\bibitem{IvanenkoETAL2020}
Y.~Ivanenko, M.~Nedic, M.~Gustafsson, B.~L.~G. Jonsson, A.~Luger, and
  S.~Nordebo, \emph{{Q}uasi-{H}erglotz functions and convex optimization},
  Royal Soc. Open Sci. \textbf{7} (2020), 191541.

\bibitem{KacKrein1974}
I.~S. Kac and M.~G. Kre\u{\i}n, \emph{R-functions--analytic functions mapping
  the upper half-plane into itself}, Amer. Math. Soc. Transl. \textbf{103}
  (1974), no.~2, 1--18.

\bibitem{Koosis1998}
P.~Koosis, \emph{Introduction to {$H_p$} spaces}, second ed., Cambridge Tracts
  in Mathematics, vol. 115, Cambridge University Press, Cambridge, 1998, With
  two appendices by V. P. Havin.

\bibitem{LugerNedic2017}
A.~Luger and M.~Nedic, \emph{A characterization of {H}erglotz-{N}evanlinna
  functions in two variables via integral representations}, Ark. Mat.
  \textbf{55} (2017), no.~1, 199--216.

\bibitem{LugerNedic2019}
\bysame, \emph{{H}erglotz-{N}evanlinna functions in several variables}, J.
  Math. Anal. Appl. \textbf{472} (2019), 1189--1219.

\bibitem{LugerNedic2019b}
\bysame, \emph{On quasi-{H}erglotz functions in one variable}, arXiv:
  1909.10198, 2019.

\bibitem{LugerNedic2020}
\bysame, \emph{Geometric properties of measures related to holomorphic
  functions having positive imaginary or real part}, J. Geom. Anal. (2020),
  DOI: 10.1007/s12220-020-00368-4.

\bibitem{Nevanlinna1922}
R.~Nevanlinna, \emph{Asymptotische {E}ntwicklungen beschr\"{a}nkter
  {F}unktionen und das {S}tieltjessche {M}omentenproblem ({G}erman)}, Ann.
  Acad. Sci. Fenn. (A) \textbf{18} (1922), no.~5, 1--53.

\bibitem{Simon1998}
B.~Simon, \emph{The classical moment problem as a self-adjoint finite
  difference operator}, Adv. Math. \textbf{137} (1998), no.~1, 82--203.

\bibitem{Vladimirov1969}
V.~S. Vladimirov, \emph{Holomorphic functions with non-negative imaginary part
  in a tubular region over a cone ({R}ussian)}, Mat. Sb. (N.S.) \textbf{79}
  (1969), 128--152, This article has appeared in an English translation [Math.
  USSR-Sb. 8 (1969), 125–146].

\bibitem{Vladimirov1979}
\bysame, \emph{Generalized functions in mathematical physics}, ``Mir'', Moscow,
  1979, Translated from the second Russian edition by G. Yankovski\u\i.

\end{thebibliography}

\end{document}